\documentclass[conference]{IEEEtran}
\IEEEoverridecommandlockouts

\usepackage{amsmath,amssymb,amsfonts}
\usepackage{algorithmic}
\usepackage{graphicx}
\usepackage{textcomp}
\usepackage{xcolor}
\usepackage{hyperref}
\usepackage{subcaption}
\usepackage{float}
\def\BibTeX{{\rm B\kern-.05em{\sc i\kern-.025em b}\kern-.08em
    T\kern-.1667em\lower.7ex\hbox{E}\kern-.125emX}}

\usepackage[style=ieee, backend=biber, giveninits=true, maxnames=5]{biblatex}
\AtEveryBibitem{%
    \clearfield{month}%
    \clearfield{language}
    \clearfield{urlyear}%
    \clearfield{issn}%
    \clearfield{doi}
    \clearfield{url}
}

\defbibenvironment{bibliography}
  {\list
     {\printfield[labelnumberwidth]{labelnumber}}
     {%
       \setlength{\labelwidth}{2em}
       \setlength{\labelsep}{0.7em}  
       \setlength{\leftmargin}{\dimexpr\labelwidth+\labelsep\relax}
       \setlength{\itemsep}{0.1\baselineskip}
       \setlength{\parsep}{0pt}
     }}
  {\endlist}
  {\item}

\DeclareBibliographyDriver{inproceedings}{%
  \printnames{author}%
  \setunit{\labelnamepunct}\newblock
  \printfield[title]{title}%
  \newunit
  \printfield{booktitle}%
  \iffieldundef{volume}{}{\setunit{\addcomma\space}\printfield{volume}}%
  \iffieldundef{pages}{}{\setunit{\addcomma\space}\printfield{pages}}%
  \setunit{\addcomma\space}\printfield{year}%
  \newunit\newblock
}

\DeclareBibliographyDriver{techreport}{%
  \printnames{author}%
  \setunit{\labelnamepunct}\newblock
  \printfield[title]{title}%
  \newunit
  \printfield{institution}%
  \iffieldundef{number}{}{\setunit{\addcomma\space}\printfield{number}}%
  \iffieldundef{pages}{}{\setunit{\addcomma\space}\printfield{pages}}%
  \setunit{\addcomma\space}\printfield{year}%
  \newunit\newblock
}

\begin{document}

\title{Modelling the Correlation Structure of Uncertain Input Parameters for Energy System Optimization\\
{\footnotesize}
\thanks{This work was partly completed within the transpAIrent.energy project, which is carried out under the Federal Ministry of Innovation, Mobility and Infrastructure Republic of Austria (BMIMI)'s 2023 “AI for Green” call for proposals. The processing is carried out by the Austrian Research Promotion Agency (FFG) on behalf of the BMIMI.\\ \\ 979-8-3195-3554-2/26/\$31.00 ©2026 IEEE}
}

\author{
\IEEEauthorblockN{
Klara Maggauer\textsuperscript{1},
Yannick Werner\textsuperscript{2},
Stefan Strömer\textsuperscript{3},
Sonja Wogrin\textsuperscript{4}
}
\IEEEauthorblockA{
\textsuperscript{1,3}Center for Energy, AIT Austrian Institute of Technology GmbH, Vienna, Austria\\
\textsuperscript{1,2,4}Institute of Electricity Economics and Energy Innovation, Graz University of Technology, Graz, Austria\\
\textsuperscript{3}Faculty of Technology, Policy and Management, Delft University of Technology, Delft, The Netherlands\\
ORCID:
\textsuperscript{1}0000-0002-5994-3201,
\textsuperscript{2}0000-0002-6674-805X,
\textsuperscript{3}0000-0002-5330-3318,
\textsuperscript{4}0000-0002-3889-7197
}
}

\maketitle

\begin{abstract}
Statistical dependence among uncertain input parameters in stochastic energy system optimization models is often ignored, even though this can substantially bias outcomes. To address this gap, we are developing a comprehensive framework for characterizing, modelling, and benchmarking statistical dependence. In this work, we present a copula-based workflow to identify, characterize, and model linear and monotonic correlation structures between input parameters, representing the first development step towards this framework. We demonstrate our workflow using solar generation, day-ahead electricity prices, and electricity demand data in Austria between 2019 and 2025. Our results show substantial linear and monotonic dependence between these variables, and that this dependence is well captured by the copula-based approach. Finally, building on this scalable foundation, we highlight key levers for next steps towards the full framework, including time-dependent and higher-dimensional dependence modelling.
\end{abstract}


\section{Introduction}\label{sec_introduction}

In energy system modelling and optimization, we often deal with predominantly statistically dependent and uncertain input parameters. Statistical dependence can, e.g., refer to correlations among integer variables, lag structures, or autocorrelation, and is due to diverse causal interrelations. Examples are an increasing share of variable renewable generation, the design of European electricity markets, demand patterns, regulatory influences, and weather- and environment-driven supply dynamics. When consistent historical datasets are available as direct inputs to the model, dependencies are inherently encoded in the data. However, for scenario generation, projections under changing system structures (e.g., increasing renewable share or evolving regulatory conditions) or when historical datasets are incomplete, it is crucial to carefully assess, quantify, and model these interrelations. Despite that, several reviews \cite{haugen_representation_2023, jordehi_how_2018, hasan_existing_2019, plaga_methods_2023} have concluded that in the majority of existing research studies applying energy system optimization models (ESOMs), this aspect is disregarded. This leads to inaccurate uncertainty representation in ESOMs and, consequently, poses a risk of misleading optimization results and poor decision-making in practice.
From a methodological perspective, a wide range of uncertainty modelling and scenario generation approaches exists for energy applications. Yet, a recurring bottleneck remains the realistic representation of joint stochastic behaviour in the presence of multiple interacting inputs \cite{roald_power_2023}. In the specific context of market and operational planning models for renewable power systems, scenario generation is repeatedly highlighted as challenging, inter alia, because multiple inputs can be statistically dependent \cite{haugen_representation_2023}. While dependence-aware approaches exist, they are not standard in ESOM workflows, and the independence assumption remains common \cite{yue_review_2018}. Copula-based methods offer a practical way to model dependence separately from the marginals and have been increasingly used in energy applications. For example, Durante et al. \cite{durante_multivariate_2022} analyse multivariate dependence between electricity prices, demand, and renewable generation using copulas. For scenario generation, several studies explicitly use copula constructions to preserve dependence structures \cite{krishna_time-coupled_2023, becker_generation_2018}. Beyond that, copula approaches have also been used to model joint price and volumetric (generation) risk in trading contexts, again emphasizing that dependence matters for realistic joint outcomes \cite{pircalabu_joint_2017}. Despite this progress, these methods are still rarely integrated into ESOM input pipelines, and straightforward comparisons against independent sampling, as well as benchmarks regarding their impact on optimization results, are often not reported.

\begin{figure}
\centering
\includegraphics[width=0.8\linewidth]{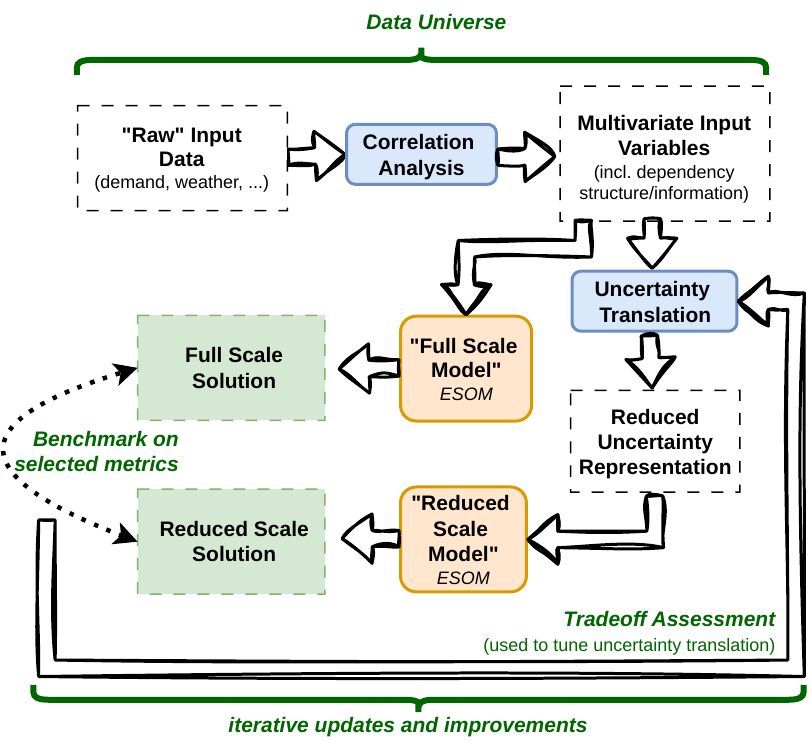}
\caption{Flow diagram of the planned full framework.}
\label{fig_int_full_framework}
\end{figure}

This gap motivates the development of a framework for the characterization, modelling, and benchmarking of statistical dependence in uncertain ESOM input parameters, as shown in Fig.~\ref{fig_int_full_framework}. We will implement this full framework step-by-step over the coming months. This paper presents a first step toward the full framework. It aims to (i) showcase the presence of linear and monotonic correlation structures between variables that are typically ESOM inputs, (ii) illustrate how these dependence structures can be modelled using copulas, and (iii) highlight key levers for future improvements of the workflow that will be addressed as the next steps towards the full framework. Specifically, in this work, we present the first step in the workflow using a representative triplet of ESOM input variables, namely solar photovoltaic (PV) generation, day-ahead electricity prices, and electricity demand in Austria between 2019 and 2025. While we do not model temporal dependence and acknowledge that in many market-clearing ESOMs prices are endogenous, this triplet is still directly relevant for settings where prices are exogenous input parameters, e.g., day-ahead bidding or self-scheduling \cite{herding_stochastic_2023}, price-based demand response \cite{al_essa_review_2025}, or price-taking dispatch problems \cite{mercier_value_2023}. We first identify and characterize the variables' empirical dependence structure using linear (Pearson) and rank-based correlation measures for monotonic and nonlinear relationships (Spearman) \cite{shemyakin_introduction_2017}. Based on this characterization, we then model their contemporaneous dependence structure using a copula-based approach, leveraging the flexibility of Gaussian and Student-t copulas to represent linear, nonlinear, and tail dependence patterns \cite{patton_review_2012, joe_dependence_2014}. Lastly, to investigate the influence of the dependencies, we pair the copulas with simple marginal models and generate joint scenarios under two settings: (i) independent sampling from the marginal distributions and (ii) dependence-preserving sampling from the copulas. We then compare these samples using the aforementioned linear and monotonic correlation measures to quantify how well the copulas preserve empirical dependence relevant to uncertain ESOM input parameters. 

The remainder of this paper is structured as follows: In Section~\ref{sec_method}, we introduce our method. We then present and discuss our findings in Section~\ref{sec_results}. Section~\ref{sec_conclusion} draws conclusions and highlights future research directions.

\section{Method}\label{sec_method}

In the following, we first present methods for characterizing the empirical dependence structure (Section~\ref{subsec_method_dependence}), estimating the marginal distributions of the three uncertain input parameters (Section~\ref{subsec_method_distributions}), and modelling their contemporaneous dependence structure (Section~\ref{subsec_method_dep_modeling}). We then describe the sampling approach in Section~\ref{subsec_method_sampling}. Lastly, the dataset used in this study and its pre-processing (Section~\ref{subsec_method_data}) are introduced. The implementation of the method is conducted in Python 3.12 \cite{python312}, leveraging packages such as \texttt{pandas} \cite{reback2020pandas}, \texttt{scipy} \cite{virtanen2020scipy}, and \texttt{statsmodels} \cite{seabold2010statsmodels}. 

\subsection{Dependence characterisation}\label{subsec_method_dependence}

To characterize the three parameters' (electricity price, load, and PV generation) correlation with each other, the following correlation measures are applied at lag zero \cite{shemyakin_introduction_2017}:
\begin{align}
    \rho_{Pearson}(X, Y) &= \frac{\text{Cov}(X,Y)}{\sigma_X \sigma_Y}, \label{eq_pearson} \\
    \rho_{Spearman}[\text{R}[X],\text{R}[Y]] &= \frac{\text{Cov}[\text{R}[X],\text{R}[Y]]}{\sigma_{\text{R}[X]} \sigma_{\text{R}[Y]}}, \label{eq_spearman}
\end{align}
where $X$ and $Y$ are the variables between which the correlation is measured, $\sigma$ is each variable's standard deviation, and $\text{Cov}$ the covariance between the variables. Ranks $\text{R}[\cdot]$ in ascending order are assigned to each observation (with ties handled appropriately). The Pearson coefficient, $\rho_{Pearson}$ \eqref{eq_pearson}, assesses linear relationships, and the rank-based Spearman coefficient, $\rho_{Spearman}$ \eqref{eq_spearman}, monotonic and nonlinear relationships. In fact, Spearman is the Pearson correlation of the paired ranks $(\text{R}[X],\text{R}[Y])$.
For both coefficients, a value close to $+1$ indicates a strong positive correlation, whereas a value near $-1$ implies a strong negative correlation. Zero correlation indicates no linear (Pearson) or no monotone (Spearman) association~\cite{shemyakin_introduction_2017}.

\subsection{Marginal distribution estimation}\label{subsec_method_distributions}

To estimate each variable's marginal distribution from the data, we collect and sort all values in ascending order for each variable $X_i$ (price, load, PV generation). From that, we extract the empirical cumulative distribution function (CDF),~$\hat{F_i}$. In general, the CDF, i.e., $F_i(x) = Pr(X_i \leq x)$, describes the probability that a parameter $X_i$ takes a value less than or equal to $x$. Using the probability integral transform (PIT) \cite{joe_dependence_2014}, we map each observation to the unit interval via $U_i=\hat{F}_i(X_i)$, enabling a representation of all variables in the common standard space $\mathcal{U}(0,1)$. Additionally, we construct the empirical inverse CDFs, $\hat{F}_i^{-1}$, as quantile functions via interpolation on the sorted values, enabling inverse transform sampling and back-transformation from uniform to physical units by $X_i=\hat{F}_i^{-1}(U_i)$.

\subsection{Dependence modelling}\label{subsec_method_dep_modeling}

Statistically dependent variables should be modelled via their joint (multivariate) distribution. Copulas are functions that describe the joint behaviour of multiple variables and separate marginal behaviour from dependence by representing the joint distribution via \cite{joe_dependence_2014}:
\begin{equation}
    F(x_1,x_2,x_3) = C(F_1(x_1),F_2(x_2),F_3(x_3)),
\end{equation}
where $F$ is the joint CDF of the three-dimensional random vector, $x_1,x_2,x_3$ are realizations of the three variables, $F_i$ are the marginal CDFs, and $C$ is the copula. In this work, marginals are handled empirically (see Section~\ref{subsec_method_distributions}), and copulas are fitted in copula space. In general, there are many copulas to choose from, each of which can capture different statistical characteristics of the joint distribution \cite{patton_chapter_2013}. Here, we model the contemporaneous dependence between the three variables using two elliptical copulas, namely the Gaussian and Student-t copulas, which provide a simple, low-parameter dependence model based on a correlation matrix. The Student-t copula additionally allows for tail dependence through its degrees-of-freedom parameter $\nu$ \cite{joe_dependence_2014}. Elliptical copulas capture central dependence patterns well but may be limited when the empirical dependence is strongly asymmetric or regime-driven \cite{joe_dependence_2014}. The copula models are fitted to a cleaned, contemporaneous dataset $X\in\mathbb{R}^{n\times 3}$, where each row corresponds to one time step and contains aligned observations of load, PV generation, and day-ahead price. Depending on the experiment, $X$ is either constructed from all hours or from a conditioned subset (only including values from 08:00 to 18:00). We then transform each variable independently to the unit interval using rank-based pseudo-observations
\begin{equation}
    U_{j,i} = \frac{\mathrm{R}(X_{j,i}) - 0.5}{n}, \qquad U_{j,i}\in(0,1),
\end{equation}
where $X_{j,i}$ denotes the observed value of variable $i\in\{1,2,3\}$, $j\in\{1,\dots,n\}$ indexes observations, $U_{j,i}$ is the corresponding pseudo-observation in copula space, and $\mathrm{R}(X_{j,i}$) is the rank of $X_{j,i}$ among all $n$ observations of variable $i$. Ranks are computed with average tie handling. To avoid numerical issues in subsequent transformations, values are clipped to $[\varepsilon,1-\varepsilon]$ with $\varepsilon=10^{-6}$. The Gaussian copula is parameterized by a correlation matrix $\Sigma$. Pseudo-observations are mapped to latent normal scores via $Z=\Phi^{-1}(U)$, where $\Phi$ denotes the standard normal CDF. The copula correlation matrix is then estimated from the pseudo-observations using a moment-based fitting routine (implemented as \texttt{fit\_corr\_param} in \texttt{statsmodels}), yielding $\hat{\Sigma}_\mathrm{Gauss}$. Given $\hat{\Sigma}_\mathrm{Gauss}$, the Gaussian copula log-likelihood is evaluated as $\sum_{j=1}^{n}\log c_\mathrm{Gauss}(U_j;\hat{\Sigma}_\mathrm{Gauss})$, whereby $c_\mathrm{Gauss}$ is the Gaussian copula probability density function (PDF). The Student-t copula is parameterized by a correlation matrix $\Sigma$ and degrees of freedom $\nu$, where a smaller $\nu$ implies stronger tail dependence. We estimate~$\nu$ by grid search over a predefined candidate set (here $\nu\in\{2,3,4,5,6,8,10,12,15,20,30,50\}$). For each candidate $\nu$, the corresponding copula correlation matrix $\hat{\Sigma}_\mathrm{Student-t}(\nu)$ is estimated using the same dependence-fitting routine \texttt{fit\_corr\_param}. We then compute the copula log-likelihood $\sum_{j=1}^{n}\log c_\mathrm{Student-t}(U_j;\hat{\Sigma}_\mathrm{Student-t}(\nu),\nu)$, where $c_\mathrm{Student-t}$ is the Student-t copula PDF, and select the $\nu$ maximizing this value. The resulting estimates are denoted $\hat{\nu}$ and $\hat{\Sigma}_\mathrm{Student-t}=\hat{\Sigma}_\mathrm{Student-t}(\hat{\nu})$. \cite{joe_dependence_2014, shemyakin_introduction_2017} Both copula fits are performed year-wise to obtain stable parameters and assess changes in dependence across years. The copula fits are discussed in the Appendix. 

\subsection{Scenario sampling}\label{subsec_method_sampling}

To isolate the effect of dependence, we generate two scenario variants with identical marginal distributions but different dependence structure: (i) independent samples and (ii) dependence-preserving samples drawn from the fitted copulas. In each case, $n=50\,000$ samples are generated. Pearson and Spearman correlation coefficients are computed on the generated samples and compared to the corresponding empirical correlations from the original dataset to assess how well linear and monotonic dependence is preserved. Sampling is performed via inverse transform sampling. First, we generate samples in the uniform space and then map them back to physical units using the empirical inverse CDFs, as explained in Section \ref{subsec_method_distributions}. For the independent scenarios, the components of $U=(U_1,U_2,U_3)$ are sampled independently with $U_i\sim\mathcal{U}(0,1)$. For the dependence-preserving scenarios, $U$ is sampled from the fitted copula (Gaussian or Student-t), yielding statistically dependent uniform samples with uniform marginals. Samples are generated i.i.d. across time steps because modelling temporal dependence is out of scope in this work.

\subsection{Data}\label{subsec_method_data}

The data used to demonstrate our workflow are the Austrian day-ahead electricity prices (EPEX market clearing prices), electricity demand, and solar PV generation from January 1, 2019, to December 31, 2025, in hourly resolution in the local timezone (CET/CEST). The data sources are the ENTSO-E transparency platform \cite{entsoe_transparency_platform} for the 2019 to 2022 data and the Austrian transmission system operator APG's transparency platform "APG Markttransparenz" \cite{apg_transparenz} for the 2023 to 2025 data. To ensure comparability across years, we divide the data into seven slices, each lasting 52 weeks, starting on a Monday and ending on a Sunday. As a result, they do not exactly correspond to the calendar years 2019 to 2025, even though, for simplicity, they are referred to as such in this work. This data is referred to as the \textit{full dataset} in the following. Additionally, we construct a \textit{conditioned dataset} by filtering the data to values between 08:00 and 18:00. This conditioning is motivated by the analysis of the PV-price correlation, which revealed that, due to the unavailability of PV during the night, using the full data biases the PV-price correlation analysis. All analyses are conducted year-wise on the seven slices.

\section{Results and Discussion}\label{sec_results}

In the following, we present and discuss our results on the characterization of dependence structures (Section~\ref{subsec_res_dep}) and the comparison between independent and dependence-preserving sampling (Section~\ref{subsec_res_sampling}) of the uncertain input parameters using the data described before.

\subsection{Dependence structure}\label{subsec_res_dep}

\begin{figure}
\centering
\includegraphics[width=\linewidth]{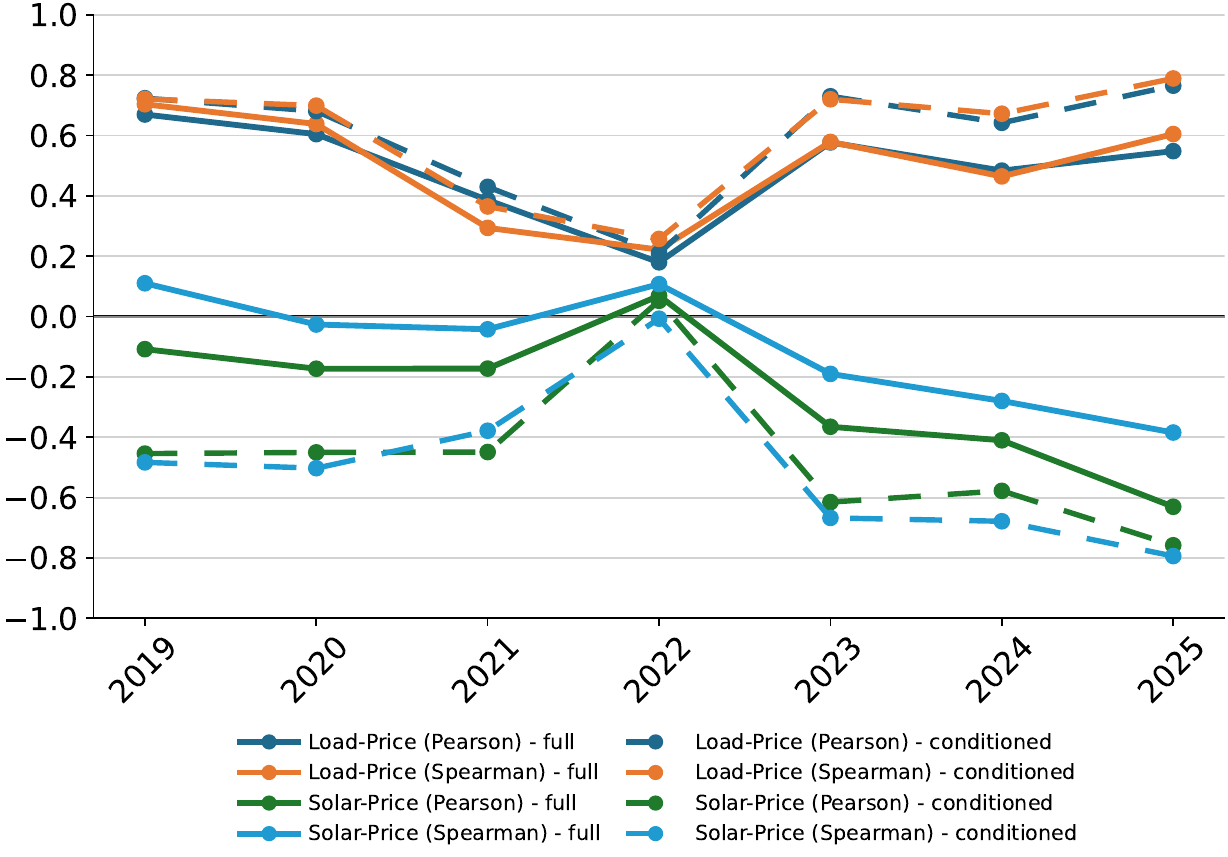}
\caption{Empirical correlation between the variables.}
\label{fig_res_emp_corr_comp}
\end{figure}

Fig.~\ref{fig_res_emp_corr_comp} illustrates the Pearson and Spearman correlation for the full (solid lines) and conditioned (dashed lines) data sets by year.

The load-price relationship shows a consistently positive correlation, but is strongly dependent on the year: In 2019/2020, it is very strong (Pearson 0.67/0.60; Spearman 0.70/0.64). This indicates a high load - high price relationship, consistent with the expectation that in high load hours, expensive fossil generators are price-setting. In 2021, the observed correlation is significantly weaker (0.39/0.29) and overall weakest (0.18/0.22) in 2022, indicating special circumstances due to the energy crisis following the Russian invasion of Ukraine. From 2023 to 2025, the correlation strengthens again (0.58/0.55). The conditioned dataset shows the same trends, with slightly more pronounced magnitudes.

Solar generation and price show a trend towards increasingly negative correlation from 2023 onwards. Here, the importance of accounting for the variable availability of PV generation throughout the day becomes obvious: Compared to the full dataset (2019-2021: Pearson ca. -0.11 to -0.17; Spearman close to zero), the conditioned dataset shows distinctly negative correlations already in 2019-2021 (Pearson/Spearman: -0.37 to -0.50). In 2022, the correlation is slightly positive but close to zero in both datasets, again suggesting special circumstances arising from the energy crisis. From 2023 to 2025, the correlation is significantly negative and increasingly strong, and again more pronounced in the conditioned (-0.58 to -0.79) than in the full (-0.19 to -0.63) dataset, consistent with increasing PV penetration and greater impact on the merit order. Finally, across the full dataset, Spearman is consistently weaker (less negative) than Pearson, indicating the effect of night hours when PV power production is zero. However, prices may still be low due to other renewable energy sources (RES) such as wind. Overall, the conditioned dataset provides better insight into the PV-price correlation.

\subsection{Independent vs. dependence-preserving sampling}\label{subsec_res_sampling}

\begin{figure*}
\includegraphics[width=0.9\textwidth]{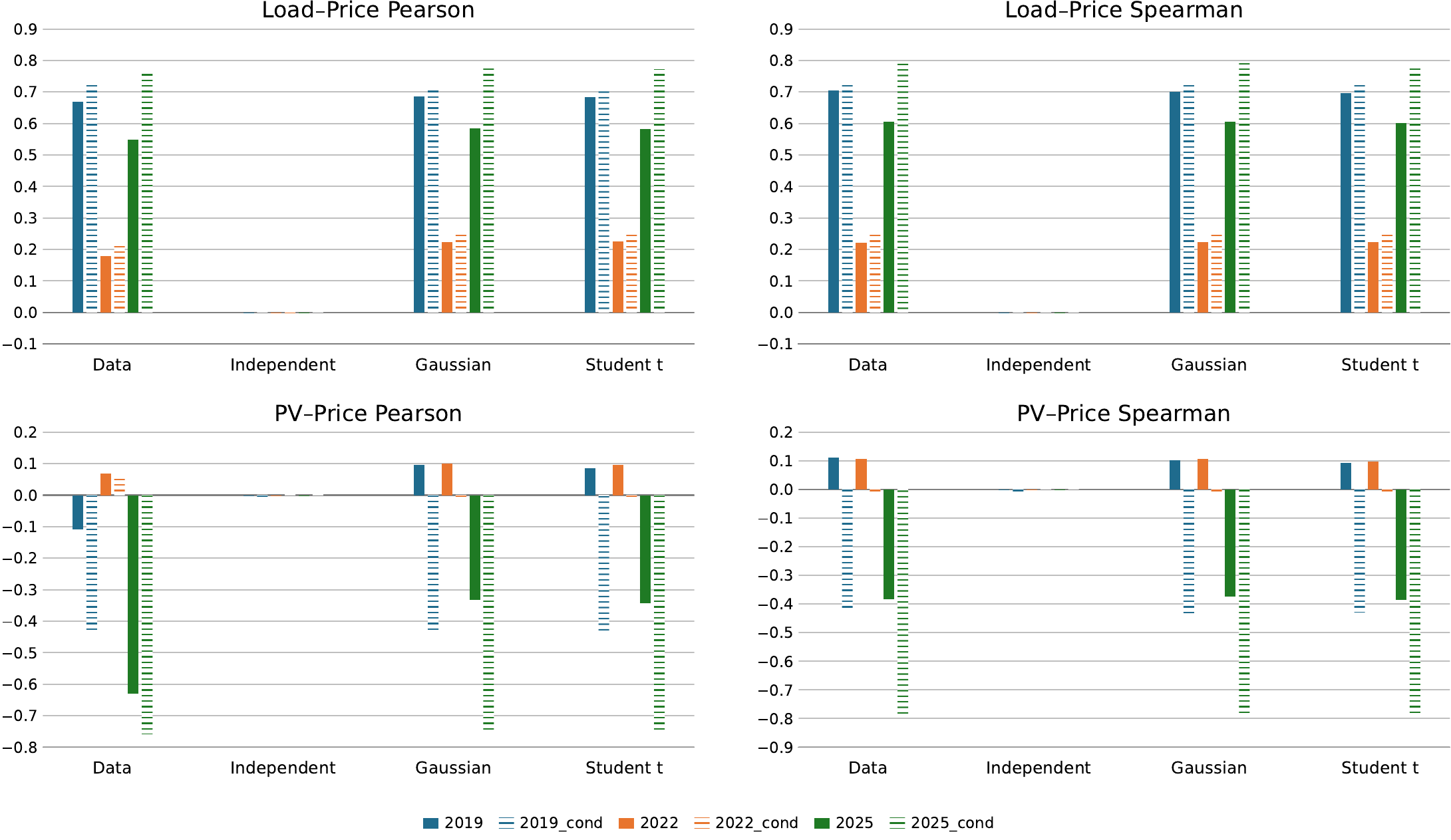}
\caption{Comparison of correlation between load-price and PV-price in the data, independent sampling, and copulas.}
\label{res_sampling}
\end{figure*}

The results of the sampling comparison measured via the Pearson and Spearman correlation coefficients, respectively, are shown in Fig.~\ref{res_sampling} for load-price and for PV-price. We select three distinct years (2019, 2022, 2025) for illustration purposes. The figures show the results for the full dataset (solid bars) as well as the conditioned dataset (striped bars) and include the empirical correlation benchmark (denoted "Data") for reference.

For the load-price case, the first observation is that independent scenarios yield approximately zero correlation, as expected. Moreover, both copulas (Gaussian and Student-t) reproduce the positive load-price relationship in the data well (Spearman coefficients are slightly closer) for both the full and conditioned datasets. Moreover, the Gaussian copula values closely match the Student-t copula values, indicating low tail dependency. The conditioned data case shows a stronger positive correlation than the full data case, with the difference becoming more pronounced in 2025.

Regarding the PV-price relationship, the copulas mostly preserve statistical dependence, whereas independent sampling yields zero correlation, as expected again. In this case, the copulas represent the linear Pearson correlation in the data significantly better when using the conditioned than the full dataset. In the latter case, in 2019, the copulas show a weak positive correlation while the empirical Pearson correlation is slightly negative. In 2025, the copulas qualitatively reflect the stronger negative correlation in the full dataset but underestimate its amplitude. This provides strong evidence of nonlinearity and seasonality, which elliptical copulas can only partially capture. However, the distinct negative correlation observed in the conditioned 2019 and 2025 datasets is well preserved by the copulas. 

Concerning the monotonic Spearman correlation, the copulas perform better even on the full dataset. This behaviour is expected because copulas are fitted on pseudo-observations in the uniform space (via the rank/PIT transform). Hence, the fitted model primarily targets the dependence structure between quantiles rather than linear co-movement in physical units. Spearman correlation is a rank-based measure and is therefore more directly tied to the copula fit, whereas Pearson correlation is sensitive to nonlinear transformations back to the original scale by the empirical inverse CDF and to intraday variability (e.g., many night hours when PV generation is zero). As a result, the copula samples can preserve the monotonic (rank) dependence well even when the linear Pearson correlation in the original units is underestimated in the unconditioned dataset. Moreover, in all cases, there is again hardly any difference in Pearson/Spearman between the Gaussian and Student-t copulas, a sign of low tail dependence.

\section{Conclusions and Outlook}\label{sec_conclusion}

In this work, we characterize the linear and monotonic dependence structure between selected uncertain ESOM input parameters (day-ahead prices, load, PV generation), model it using Gaussian and Student-t copulas, and benchmark copula-based sampling against independent sampling.

In conclusion, our results show that linear and monotonic dependencies exist between the considered parameters. For accurate representation in ESOMs, these dependencies should therefore be explicitly modelled, e.g., using copula-based methods, which we have shown to maintain dependence structures well. The presented method can be readily extended to incorporate additional parameters, such as other electricity generation technologies from intermittent RES.
Furthermore, the observed dependence structure is clearly year-dependent: The PV-price case, in particular, becomes increasingly negative from 2023 onwards due to increasing PV penetration and the resulting effect on the merit order. In addition, the energy crisis is recognizable as a special case. These findings imply that system changes, such as variations in RES generation shares or extreme events, strongly affect correlation structures and thus should be considered in statistical dependence analyses.
When considering the full dataset, the elliptical copulas applied in this work capture the central trends of the investigated parameters, but not all details, especially for PV-price. This highlights the limits of simple elliptical copulas for data with strong structural patterns. Conditioning the dataset to only contain values from 08:00 to 18:00 improves copula performance for PV-price, especially with regard to preserving monotonic Spearman correlation. This indicates nonlinearity and seasonality in the investigated parameters, which should be addressed. Further, when comparing the results for the Gaussian and Student-t copulas, hardly any differences are discernible in their shapes (high degrees of freedom value $\nu$). Therefore, there is no evidence for strong tail dependence, at least as measured by the calculated linear and monotonic correlations.
Moreover, our work shows that independent sampling effectively eliminates statistical dependencies. In contrast, samples generated with Gaussian and Student-t copulas preserve linear and monotonic correlations well. This suggests that even simple elliptical copulas are a powerful tool for preserving the dependence structure.\\


In future work, the most crucial addition is to quantify the practical implications for actual ESOM results by embedding the generated scenarios in an optimization model and comparing outcomes when dependencies are represented or ignored. This would extend the present, deliberately simple benchmark against independent sampling.
Moreover, temporal dependence in the considered time series data should be modelled, as so far, only the simultaneous dependence structure has been considered, neglecting auto- and cross-correlations. This can be done using time-varying copulas, e.g., GARCH (Generalized Autoregressive Conditional Heteroskedasticity)-based copula models or other types of dynamic copulas \cite{patton_chapter_2013}.
In addition, seasonality needs to be addressed more comprehensively: Our results of the conditioned vs. the full dataset show that dependencies are influenced by daily patterns, particularly in the case of PV generation. It is expected that there are also weekly/seasonal effects, which again emphasizes the importance of explicitly modelling the temporal dependence structure in future work. In this context, further conditioning or focusing on residualized (de-seasonalized) datasets \cite{box_box_2013} should be done.
Furthermore, nonlinearity, regime effects, and tail dependence need to be investigated and modelled in greater detail: So far, the focus has been on elliptical copulas, which capture central trends but fail to represent asymmetric effects. To address this, possible extensions include vine copulas, regime-switching copulas, or copulas per time window/season, combined with comprehensive goodness-of-fit testing, e.g., via Cramer-von Mises or Rosenblatt transforms \cite{patton_chapter_2013}. In that vein, future work should also include a comprehensive analysis of changes in correlations during extreme events \cite{shemyakin_introduction_2017}.
Another aspect to be addressed in future work is the robustness of the presented analyses over different years with varying system states. This could be achieved, for example, by fitting separately normal years and the crisis year 2022, and by varying the RES expansion level.

\printbibliography

\section*{Appendix}\label{sec_appendix}

\begin{figure}[H]
    \centering

    \begin{subfigure}{0.4\linewidth}
        \centering
        \includegraphics[width=\linewidth]{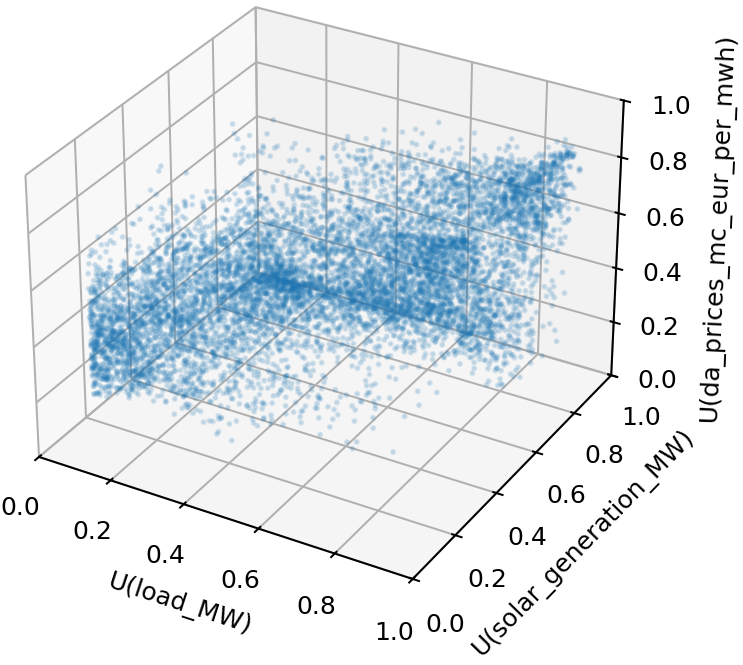}
        \caption{Empirical distribution.}
    \end{subfigure}
    \hfill
    \begin{subfigure}{0.4\linewidth}
        \centering
        \includegraphics[width=\linewidth]{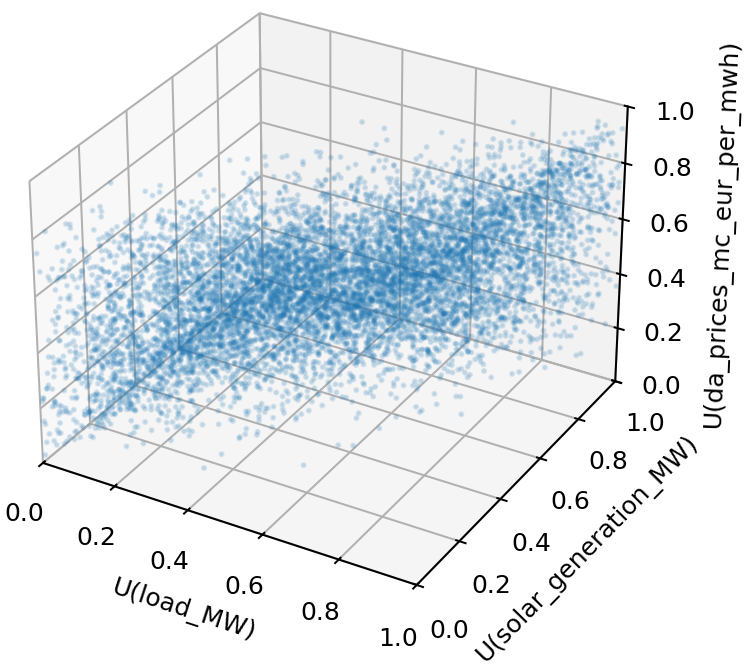}
        \caption{Gaussian copula.}
    \end{subfigure}

    \vspace{0.1em}

    \begin{subfigure}{0.4\linewidth}
        \centering
        \includegraphics[width=\linewidth]{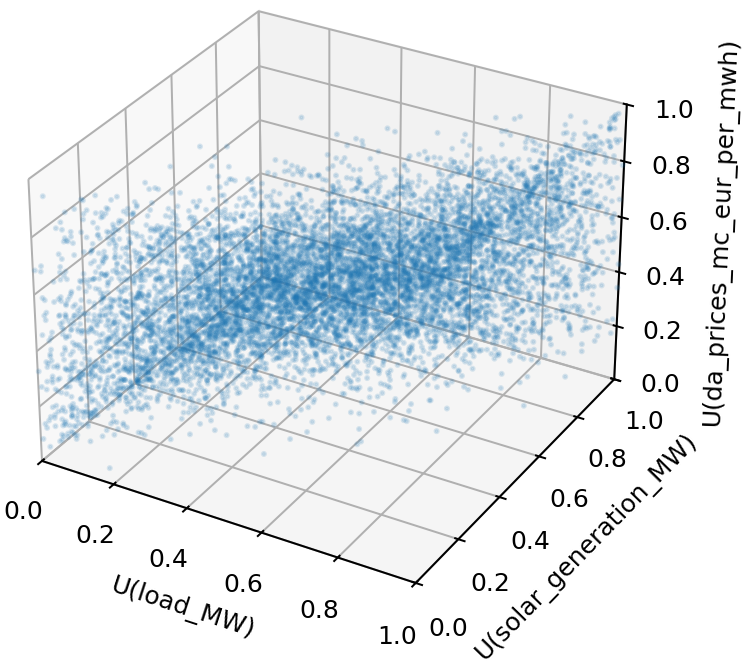}
        \caption{Student-t copula.}
    \end{subfigure}

    \caption{Samples drawn from slice 7 (2025) of the full dataset in U-space.}
    \label{fig_copulas}
\end{figure}

Fig. \ref{fig_copulas} illustrates the empirical multivariate distribution and samples drawn from the copulas in U-space. In the empirical U-space, the point cloud does not resemble a simple elliptical shape but rather forms a structured sheet. This indicates that the joint dependence is not governed by a single homogeneous regime, but is strongly shaped by daily, weekly, and seasonal patterns (e.g., day vs. night, weekday vs. weekend, summer vs. winter). In particular, PV generation is exactly zero or near zero for many hours (nighttime and low-irradiance periods). When transformed to U-space, this creates a high concentration of points in a limited solar-quantile region, which contributes to the geometry seen in the empirical cloud. Both fitted copulas reproduce the broad global pattern of the dependence, i.e., the overall direction of co-movement implied by the fitted correlation structure. This indicates that copula-based scenarios should preserve Pearson/Spearman correlations as opposed to independent sampling. However, the copula samples look smoother and more elliptical and therefore wash out the geometry visible in the empirical U-space. This mismatch is consistent with the fact that elliptical copulas can approximate central dependence well but have limited flexibility for asymmetric dependence. Finally, the Gaussian and Student-t copula samples appear very similar in this case. This suggests that, for this dataset and year, the additional tail-dependence flexibility of the Student-t-copula is either weakly supported (e.g., the fitted $\nu$ are relatively large, so the Student-t-copula behaves close to Gaussian) or not strongly expressed in the empirical U-space.

\end{document}